\documentclass[12pt]{amsart}

\usepackage[english]{babel}

\usepackage{amsmath}

\usepackage{changes}

\usepackage{amsthm}

\usepackage{amsfonts} 

\usepackage{hyperref}

\usepackage{epsfig}

\usepackage{amssymb}

\usepackage{environ}

\usepackage{mathrsfs}

\usepackage{caption}

\usepackage{tikzscale}

\usepackage{tikz-cd}

\usetikzlibrary{cd}

\usetikzlibrary{patterns}

\usepackage{enumitem}

\usepackage{dsfont}

\usepackage{aligned-overset}

\usepackage{mathtools}

\usepackage{float}

\usepackage{subcaption}

\usepackage{yfonts}

\usepackage[
backend=biber,
style=alphabetic,
]{biblatex}

\usepackage{import}

\tikzcdset{scale cd/.style={every label/.append style={scale=#1},
    cells={nodes={scale=#1}}}}

\DeclareMathOperator{\krn}{Ker}
\DeclareMathOperator{\im}{Im}

\newtheorem{theorem}{Theorem}[section]
\newtheorem{definition}[theorem]{Definition}
\newtheorem{corollary}[theorem]{Corollary}
\newtheorem{proposition}[theorem]{Proposition}
\newtheorem{lemma}[theorem]{Lemma}
\newtheorem{remark}[theorem]{Remark}

\theoremstyle{definition}
\newtheorem{example}[theorem]{Example}

\newcommand{\Ocal}{\mathcal{O}}

\DeclareMathAlphabet{\mathcal}{OMS}{cmsy}{m}{n}

\usepackage{fullpage}

\title{A multitor formula for a non-lci intersection}

\author{Oscar Finegan}

\begin{document}

\maketitle{}
\begin{abstract}
    In this paper we produce the first known formula for cohomologies of the derived tensor products of structure sheaves of subschemes in the case where the intersection of the subschemes is not a local complete intersection. The case covered here is where the intersection instead consists of two local complete intersection components, one of codimension 1 and the other of arbitrary codimension.
\end{abstract}

\section{Introduction}

Classically, the structure sheaf of the intersection of two subschemes is isomorphic to the tensor product of the structure sheaves of the subschemes. This enables us to study the non-reduced structures on the intersection, and distinguish between the scheme structure arising as an intersection and merely imbedding as a subscheme into our ambient space. However, we can obtain more data from the intersection, and it is pertinent geometrically. Indeed, in \cite{Serre}, Serre gave his \textit{intersection multiplicity formula;}

\begin{theorem}[\cite{Serre}]
    Let $X$ be a regular variety, and $U,V$ subvarieties of complementary dimension whose intersection is 0-dimensional. Then at a point $P \in U\cap V$, the intersection multiplicity is given by; 
        \[m_P(U,V) = \sum_{i = 0}^\infty (-1)^i length_{\Ocal_{X,P}}(Tor_{\Ocal_{X,P}}^i(\Ocal_{U,P},\Ocal_{V,P}))\]
\end{theorem}
In the local complete intersection (lci) case, the higher Tors disappear and this agrees with the classical definition of intersection multiplicity, but in the non-lci case all of the higher Tors are needed as correction terms in order for the algebraic formula to give the correct numerics for the intersection multiplicities. For example, the intersection in $\mathbb{A}^4$ of a union of two planes meeting in a point (e.g. $U = (xz, xw, yz, yw)$) with a general plane through that point (e.g. $V = (x - z,y - w)$)  should have intersection multiplicity 2, but $length(\Ocal_{U,P}\otimes \Ocal_{V,P}) = 3$, so we need $Tor^1(\Ocal_{U,P},\Ocal_{V,P}) = 1$ to correct it. 
\par
Presently, derived algebraic geometry tells us that to have good formal properties for our intersection theory, we should work with \textit{derived schemes}. These can be thought of as topological spaces with sheaves of commutative dg-algebras as structure sheaves, with ordinary schemes having their structure sheaf concentrated in degree 0 \cite{lurie2009derived}\cite{toën2014derived}\cite{https://doi.org/10.48550/arxiv.1412.5233}. The derived scheme associated to an intersection of derived schemes is their set theoretic intersection equipped with the derived tensor product of their structure sheaves. Even for an ordinary non-derived scheme $X$, the intersection of subschemes $Y_1,\dots,Y_n$ should have as its structure sheaf a dg-algebra which is isomorphic in $DCoh(X)$ to $\Ocal_{Y_1}\otimes^L \dots \otimes^L \Ocal_{Y_n}$. It is therefore important to be able to compute its cohomologies, these are the \textit{multitors} of $Y_1,\dots, Y_n;$
\[Tor^q_{\Ocal_X}(\Ocal_{Y_1},\dots,\Ocal_{Y_n}) := H^{-q}(\Ocal_{Y_1}\otimes^L\dots\otimes^L\Ocal_{Y_n}).\]
\par
In this paper, we give an explicit formula for these multitors in a case where the lci $Y_i$ intersect in a non-lci subscheme. To compute the multitors one replaces the sheaves $\Ocal_{Y_i}$ by flat resolutions. Theoretically, flat resolutions always exist, but in practice computing the multitors without explicit models is very difficult. This leads us to study the case where the $Y_i$ are local complete intersections. In this case, locally the $\Ocal_{Y_i}$ have Koszul resolutions. These are explicit resolutions by free sheaves which make it straightforward to compute the multitors. Normally, one has to provide a gluing argument to show how the local models of the Koszul cohomology glue to a global sheaf. If the intersection of the $Y_i$ is also lci, we are gluing free local models into a locally free global answer as computed by Scala in \cite{https://doi.org/10.48550/arxiv.1510.04889};
\begin{theorem}{\cite{https://doi.org/10.48550/arxiv.1510.04889}}\label{Scala}
Let $X$ be a smooth algebraic variety and $Y_1,\dots,Y_n$ be locally complete intersection subvarieties of $X$ such that the intersection $Z:= Y_1\cap\dots\cap Y_n$ is also a locally complete intersection. Then:
\[Tor_{\Ocal_X}^q(\Ocal_{Y_1},\dots,\Ocal_{Y_n}) = \bigwedge^q\mathcal{E}_Z\]
where $\mathcal{E}_Z$ is the conormal excess bundle \[\mathcal{E}_Z := \krn\left(\bigoplus (\mathcal{N}_{Y_i/X})^\vee|_Z \rightarrow \mathcal{N}_{Z/X}^\vee\right).\]
\end{theorem}
\noindent
An intersection of subschemes also has an associated  
\textit{normal excess sheaf} \[\mathcal{E}_Z^{norm}:=\left(\bigoplus^l_{i=1}\mathcal{N}_{Y_i}|_Z\right)/\mathcal{N}_Z.\]
  In the lci case both of the associated excess sheaves are locally free and they are dual to each other. The main result of this paper extends Theorem \ref{Scala} above to a case where $Z$ is not lci. In this case the dual of the conormal excess sheaf will contain the normal excess sheaf, but nothing can be said for the dual of the normal excess sheaf. Henceforth, excess sheaf shall mean conormal excess sheaf. The case we are considering is where the intersection can be stratified by codimension into two parts $Z$ and $D$, where $D$ is of codimension 1 and $Z$ is lci. Here, the $Y_i$ are then forced to be divisors so we have global Koszul complexes whose cohomologies model the multitors of the $\Ocal_{Y_i}$ and the $\Ocal_{Y_i - D}$. Denote the differential of this second Koszul complex by $\delta^q$. We prove the following;
\begin{theorem}[Main Theorem]
    Let $X$ be a non-singular variety with lci subvarieties $Y_1,\dots,Y_n$. Suppose that the intersection $\bigcap Y_i = Z\cup D$ where $Z$ and $D$ are lci with $codim(D,X) = 1$. Then 
    \[Tor^q_{\Ocal_X}(\Ocal_{Y_1},\dots, \Ocal_{Y_n}) \cong \mathcal{H}^q \otimes \Ocal(-qD)\]
    where $\mathcal{H}^q$ is the categorical pullback in $Mod_X$ in the diagram;
    \begin{center}
    \begin{tikzcd}
&\mathcal{H}^q \arrow[r]\arrow[d] &\krn(\delta^{-q})|_D \arrow[d] \\
&\bigwedge^q\mathcal{E}_Z \arrow[r] &\bigwedge^q\mathcal{E}_Z|_{D\cap Z}\\
\end{tikzcd}
\end{center}
 Here, $\mathcal{E}_Z$ is the conormal excess bundle for $Z$ as the intersection of divisors $Y_i-D$, and $\delta^q$ is the differential of the associated Koszul complex for $Z$. 
\end{theorem}
In effect, this result tells us that the multitor of the $Y_i$ is given by a gluing of excess sheaves on each component of the intersection. We make explicit this gluing, in that we give an object on the intersection of the components and maps from the excess sheaves on the components to this object over which the multitor is the pullback. On the open sets away from each component, our subschemes satisfy the hypotheses of Theorem \ref{Scala} and therefore the multitors on these open loci are excess bundles. In our case, the intersection of $Z$ in $D$ has codimension at least 2 in $D$ so these locally free sheaves have unique extensions to all of $D$, however Example \ref{Example 2} demonstrates that the naive gluing of these sheaves will not yield the correct answer and these more complicated kernel sheaves are required. This example also demonstrates that the multitors do not have an exterior algebra structure.  
\par
The structure of this paper is as follows. In section 2 we give an overview of the results we need about Koszul complexes and their cohomologies, as well as several results on local complete intersection subschemes. In section 3 we prove several algebraic results which approximate and motivate our main result, and then provide a proof for the main theorem. Section 4 provides some ideas for further directions. 

\subsection*{Acknowledgements}

The author is grateful to Timothy Logvinenko for many enlightening discussions on the subject of this paper. 

\section{Preliminaries}

Our convention throughout this paper is that rings are commutative with unit and schemes are noetherian, seperated and of finite type over a field $k$. 
\subsection{Koszul complexes}

We begin with an aside on Koszul complexes. Most of the results here can be found in \cite{Matsumura} for the algebraic case or \cite{F-L} for the geometric case. 

\subsubsection{Koszul complexes in commutative algebra}\label{K-alg}

 Let $R$ be a ring and let $s: M \rightarrow R$ be a map of $R$-modules. Then there is an induced map $1\otimes\dots\otimes1\otimes s: M^{\otimes p} \rightarrow M^{\otimes p-1}$. Now $\bigwedge^p M$ naturally embeds into $M^{\otimes p}$ given by linearly extending the mapping 
 \[e_1\wedge \dots \wedge e_p \mapsto \sum_{\sigma \in S_p}(-1)^{sgn(\sigma)}e_{\sigma(1)}\otimes \dots \otimes e_{\sigma(p)}.\]
 The image of $1\otimes \dots \otimes 1 \otimes s$ is contained in the image of the inclusion of $M^{\otimes p-1}$. We denote the induced map on the level of exterior powers by $1\wedge \dots \wedge s: \bigwedge^pM \rightarrow \bigwedge^{p-1}M.$  We sometimes abuse notation and write this map as $1\wedge s$.
 
 \begin{definition}
 The \textit{Koszul complex} of $(M,s)$ over $R$ is the complex of $R$-modules concentrated in negative degree whose terms are $K^{-p}(M,s) = \bigwedge_R^p M$ for $p \geq 0$ with differential $1 \wedge \dots \wedge 1\wedge s$ (we will identify $ \bigwedge_R^{p-1}M\wedge_R R$ with $\bigwedge_R^{p-1}M$ and drop the subscripts).
 \end{definition}
 \noindent
Explicitly, these differentials act by; 
\[m_1\wedge \dots \wedge m_p \mapsto \sum_{i =1}^p (-1)^i s(m_i) (m_1\wedge \dots \wedge\widehat {m_i} \wedge\dots \wedge m_p)\]
where $\widehat{m_i}$ means to omit the $m_i$ term from the exterior product.
Of interest to us are the cohomologies of this complex, in particular, we are interested in the case when $M$ is a free $R$-module of finite rank. We write $K^\bullet(R^r,s) =: K^\bullet(f_1,\dots,f_r) = K^\bullet(\underline{f})$ where $f_i = s(e_i)$ is the image of the $i^{th}$ standard generator. 
\par
Koszul complexes are nicely behaved with respect to tensor products;
\begin{proposition}
Suppose we have two linear maps $s: E \rightarrow R, t : F \rightarrow R$ of free modules of finite rank into $R$. Then we have the identification 
\[K^\bullet(E,s)\otimes K^\bullet(F,t) \cong K^\bullet(E\oplus F,s\oplus t).\]
In particular  if $F \cong R$, we get
\[K^\bullet(E\oplus R,s\oplus t) \cong Cone(t:K^\bullet(E,s)\rightarrow K^\bullet(E,s)).\]
\end{proposition}
The following is a well known criterion for a Koszul complex $K^\bullet(\underline{f})$ to be a resolution of $R/((f_1,\dots,f_n)$; 

\begin{definition}
We say a sequence of elements $f_1,\dots, f_n \in R$ form a regular sequence in $R$ if $f_1$ is not a zero-divisor in $R$ and for each $i > 1$, $f_i$ is not a zero-divisor in $R/(f_1,\dots,f_{i-1})R.$
\end{definition}

\begin{proposition}
Let $f_1,\dots,f_r$ form a regular sequence in a ring $R$. Then $H^i(K^\bullet(f_1,\dots,f_r)) = 0$ for all $i \leq {-11}$ and $H^0(K^\bullet(f_1,\dots,f_r)) = R/(f_1,\dots,f_r)$, i.e. $K^\bullet(f_1,\dots,f_r)$ is a free resolution of $R/(f_1,\dots,f_r).$
\end{proposition}

In general the converse does not hold (see \cite{Kabele}), but for instance in Noetherian local rings with $f_1,\dots, f_n$ in the maximal ideal it does \cite{Matsumura}. 

\subsubsection{Koszul complexes in Geometry}
The algebraic notions of \S\ref{K-alg} can be interpreted in an algebro-geometric context. 
\par
Let $X$ be a scheme, $\mathcal{E}$ an $\Ocal_X$-module and $s : \mathcal{E} \rightarrow \Ocal_X$ a map of $\Ocal_X$-modules. We can define the Koszul complex $K^\bullet(\mathcal{E},s)$ analogously to the algebraic case. It is a complex concentrated in negative degrees with $K^{-p}(\mathcal{E},s) = \bigwedge^p\mathcal{E}$ for $p\geq 0$ and differential $1 \wedge \dots \wedge 1\wedge s$. We work in the case that $\mathcal{E}$ is a locally free $\Ocal_X$ module, where locally $K^\bullet(\mathcal{E},s)$ is isomorphic to a complex of the form $K^\bullet(f_1,\dots,f_r).$

\begin{proposition}
Suppose we have two morphisms $s: \mathcal{E} \rightarrow \Ocal_X, t : \mathcal{F} \rightarrow \Ocal_X$ of locally free $\Ocal_X$-modules of finite rank into $\Ocal_X$. Then we have an identification 
\[K^\bullet(\mathcal{E}\oplus \mathcal{F},s\oplus t) \cong K^\bullet(\mathcal{E},s)\otimes K^\bullet(\mathcal{F},t).\]
Additionally, if $t: \mathcal{L} \rightarrow \Ocal_X$ is a map from a line bundle, then we have 
\[K^\bullet(\mathcal{E}\oplus \mathcal{L} ,s\oplus t) \cong Cone(t:K^\bullet(\mathcal{E},s)\otimes \mathcal{L}\rightarrow K^\bullet(\mathcal{E},s)).\]
\end{proposition}

\begin{definition}
We call a section $s: \mathcal{E}\rightarrow \Ocal_X$ of a locally free sheaf $\mathcal{E}^\vee$ regular if the image ideal sheaf $\mathcal{I}_s$ can be locally generated by a regular sequence of sections of $\Ocal_X$. We call an ideal sheaf with this property a regular ideal sheaf.
\end{definition}

\begin{proposition}
For any regular section $s: \mathcal{E} \rightarrow \Ocal_X$, $K^\bullet(\mathcal{E},s)$ is a global locally free resolution of the structure sheaf of $Z(s)$, the zero locus of the section $s$.
\end{proposition}

\noindent
Later we will need the following technical lemma on the interaction of exterior products with tensor products of locally free sheaves. 

\begin{lemma}
\label{lem:line}
    Let $\mathcal{E}, \mathcal{L}$ be locally free sheaves of ranks $n$ and $1$ respectively on a scheme $X$. Then
    \[\bigwedge^k \left(\mathcal{E}\otimes \mathcal{L}\right) \cong \left(\bigwedge^k \mathcal{E}\right) \otimes \mathcal{L}^{\otimes k}.\]
\end{lemma}

\begin{proof}
There is a natural morphism $(\mathcal{E} \otimes \mathcal{L})^{\otimes k} \rightarrow \bigwedge^k \mathcal{E} \otimes \mathcal{L}^{\otimes k}$. This descends to a morphism $\bigwedge^k(\mathcal{E}\otimes \mathcal{L}) \rightarrow \bigwedge^k \mathcal{E} \otimes \mathcal{L}^{\otimes k}$ since any section $(e \otimes l_1) \wedge (e \otimes l_2)$ is locally everywhere 0 (as $\mathcal{L}$ has rank 1) and hence globally 0. This morphism is clearly locally everywhere an isomorphism as $\mathcal{L}$ is an invertible sheaf and therefore is a global isomorphism. 
\end{proof}

\subsection{Regular immersions and local complete intersections}
 Throughout this section $X$ is a nonsingular variety, although many of the results will hold if $X$ is a Cohen-Macaulay scheme. 

\begin{definition}
    We call an immersion $i: Y \rightarrow X$ regular if the ideal sheaf $\krn(\Ocal_X \rightarrow i_*\Ocal_Y)$ is a regular ideal sheaf. 
\end{definition}

When $X$ is Cohen-Macaulay then the notion of regular immersion coincides with local complete intersection \cite{hartshorne_1977}. 
\begin{definition}
A closed subscheme $Y \subset X$ is called a local complete intersection if locally around every closed point $y \in Y$, there is a neighbourhood $U$ of $y$ in $X$ such that $\mathcal{I}_Y(U)$ is generated by codim($Y,X$) elements. 
\end{definition}

Note that this definition makes use of a global notion of codimension and not a local one. An immediate consequence is that if $Y$ is not equidimensional (hence not equicodimensional) then it cannot be a lci.  

\begin{proposition}
The ideal sheaf of a local complete intersection is locally generated by a regular sequence of length equal to the codimension. 
\end{proposition}

A local complete intersection can be seen in one way as a generalisation of a smooth subvariety, and in another as a generalisation of an effecive Cartier divisor. Indeed, every smooth subvariety and every effective Cartier divisor is a local complete intersection subvariety. Local complete intersections and smooth subvarieties are the only subvarieties whose cotangent complex has finite amplitude. Effective Cartier divisors are simply codimension 1 local complete intersections. We will use the terms regular immersion and lci interchangeably since we are working with a smooth ambient space. 
Examples of regular immersions include closed immersions of smooth varieties and effective Cartier divisors. 
\par
It is possible to have non-lci subvarieties of smooth varieties which are equidimensional. They must be singular, but there are simple examples such as; 

\begin{example}
Let $X = \mathbb{A}^3_k$ and let $Y = Z(xy,xz,yz)$ the closed subscheme given as the union of the 3 coordinate axes. Then $Y$ is not a local complete intersection in $X$. Indeed, at the origin, 3 equations are required to carve out $Y$, but $Y$ only has codimension 2. 
\end{example}

For our results we will need to know about the associated points of a local complete intersection. We imagine this result is known but include a short proof.

\begin{lemma}\label{no-embedded}
    Let $Z$ be a local complete intersection subscheme of a noetherian Cohen-Macaulay (CM) scheme $X$. Then the only associated points for $Z$ are the generic points of its irreducible components. 
\end{lemma}

\begin{proof}
    By \cite{Matsumura}[Thm 17.3], if $A$ is a local, noetherian, CM ring and $f_1,\dots,f_n$ are a regular sequence in $A$, then $A/(f_1,\dots,f_n)$ is a CM-module over $A$. Additionally, CM modules have no embedded associated primes. Now, by \cite{Matsumura}[Thm 6.2] if a module $M$ over $A$ has the property that $M_P$ has no embedded primes over $A_P$ for every prime ideal $P$, then $M$ has no embedded prime ideals over $A$. Hence the structure sheaf of $Z$ has no embedded associated points, and therefore its associated points are only the generic points of its irreducible components. 
\end{proof}



\section{Multitors}
We want to investigate the multitors of intersections which are not lci. An obvious example to consider is where the intersection is not equidimensional but in each codimension the corresponding component is lci. As the complexity of the geometry increases with the codimension, the most straightforward case to look at is where we have two components to the intersection, where both are lcis but one has codimension 1. We begin with an algebraic result.

\begin{proposition}
\label{prop: affine}
Let $R$ be a ring and $f_1,\dots,f_n$ a sequence in $R$. Let $x$ be a non-zero-divisor in $R$ and consider the sequence $xf_1,\dots,xf_n$. We calculate the Koszul cohomology of this sequence 
\[H^q(K^\bullet(xf_1,\dots,xf_n)) \cong \krn(\partial_{\underline{f}}^q)/x\cdot\im(\partial_{\underline{f}}^{q-1})\]
where $\partial_{\underline{f}}$ denotes the differential in $K^\bullet(f_1,\dots,f_n)$.
\end{proposition}

\begin{proof}
Let $\partial_{\underline{xf}}$ denote the differential in $K^\bullet(xf_1,\dots,xf_n)$. The first thing to note is that we have $\partial_{\underline{xf}} = x\partial_{\underline{f}}$ as can be seen from the explicit formulation of the differential in a Koszul complex. Then since $x$ is a non-zero-divisor, we have 
\[\krn(\partial_{\underline{xf}}^q) = \krn(\partial_{\underline{f}}^q), \quad\im(\partial_{\underline{xf}}^q) = x(\im(\partial_{\underline{f}}^q)). \]
We therefore have the identifications
\begin{align*}
    H^q(K^\bullet(xf_1,\dots,xf_n)) &= \krn(\partial_{\underline{xf}}^q)/\im(\partial_{\underline{xf}}^{q-1})\\
    &= \krn(\partial_{\underline{f}}^q)/x(\im(\partial_{\underline{f}}^{q-1}))\\
\end{align*}
as stated.
\end{proof}

\begin{corollary}
    If the sequence $f_1,\dots,f_n$ is regular, then for any $q<0$ 
    \[H^q(K^\bullet(xf_1,\dots,xf_n)) \cong \krn(\partial_{\underline{f}}^q) \otimes R/(x).\]
\end{corollary}
\begin{proof}
    The sequence $f_1,\dots,f_n$ being regular implies that  $\im(\partial^{q-1}_{\underline{f}}) = \krn(\partial^q_{\underline{f}})$ for all $q<0$ so 
    \begin{align*}
        H^q(K^\bullet(xf_1,\dots,xf_n)) &= \krn(\partial_{\underline{f}}^q)/x(\krn(\partial^q_{\underline{f}}))\\
        &\cong \krn(\partial_{\underline{f}}^q) \otimes R/(x)
    \end{align*}
\end{proof}
We'd like to relate $\krn(\partial_{\underline{f}}^q)/x(\im(\partial_{\underline{f}}^{q-1}))$ with the cohomology of $K^\bullet(f_1,\dots,f_n)$ for which we make use of the following general algebraic statement. 

\begin{lemma}\label{II + III}
    Let $M$ be an $R$-module with submodules $P,Q$. Then the commuting square of projections 
    \begin{tikzcd}
        &M/(P\cap Q) \ar[r]\ar[d] &M/Q \ar[d] \\
        &M/P \ar[r] &M/(P+Q) 
    \end{tikzcd}
    is cartesian.
\end{lemma}
\begin{proof}
    The given square of projections commutes so there is an induced map to the pullback $M/(P\cap Q) \rightarrow X$. An element of $X$ is uniquely determined by a pair $(\pi_P(m),\pi_Q(n))$ such that $\pi_{P+Q}(m) =\pi_{P+Q}(n).$ Then $n = m + p + q$ so $(\pi_P(m),\pi_Q(n)) = (\pi_P(m+p),\pi_Q(m+p))$. Hence $M/(P\cap Q) \rightarrow X$ is surjective. Suppose that two elements $m,n$ of $M$ map to the same element of $X$ via $M/(P\cap Q)$. Then $\pi_P(m) = \pi_P(n), \pi_Q(m) = \pi_Q(n)$ so $\pi_{P \cap Q}(m) = \pi_{P \cap Q}(n)$. Hence our map $M/(P\cap Q) \rightarrow X$ is injective and hence an isomorphism.
\end{proof}

\begin{remark}
    Note that \ref{II + III} holds for sheaves of modules on a topological space as the maps are global but checking the isomorphism is local on the stalks. 
\end{remark}

\begin{proposition}\label{alg-pullback}
    In the context of Proposition \ref{prop: affine}, suppose additionally that $x$ does not belong to any prime ideal in $Ass_R(\krn(\partial_{\underline{f}}^q)/\im(\partial_{\underline{f}}^{q-1}))$. Then the object $\krn(\partial_{\underline{f}}^q)/x\cdot\im(\partial_{\underline{f}}^{q-1})$ fits into a cartesian square of projections
    \[
    \begin{tikzcd}
        &\krn(\partial_{\underline{f}}^q)/x\cdot\im(\partial_{\underline{f}}^{q-1}) \ar[r]\ar[d] &\krn(\partial_{\underline{f}}^q) \otimes R/(x) \ar[d] \\
        &\krn(\partial_{\underline{f}}^q)/\im(\partial_{\underline{f}}^{q-1}) \ar[r] & \krn(\partial_{\underline{f}}^q)/\im(\partial_{\underline{f}}^{q-1})\otimes R/(x).
    \end{tikzcd}
    \]
\end{proposition}
\begin{proof}
     For the sake of notation, denote $\krn(\partial_{\underline{f}}^q)$ and $\im(\partial_{\underline{f}}^{q-1})$ by $\krn$ and $\im$ respectively. In Lemma \ref{II + III}, take $M$ to be $\krn/x\cdot\im$, $P = \im/x\cdot\im$, $Q = x \cdot \krn/x\cdot\im$. Since $x$ does not belong to any element of $Ass_R(\krn/\im)$, we can equivalently say that $(\im) \cap (x\cdot \krn) = x\cdot \im$. To see this note that if $xy \in \im, y \in \krn$, then $\overline{xy} = 0$ in $\krn/\im$ but the action of $x$ does not kill any elements of $\krn/\im$ so $\overline{y} = 0 \in \krn/\im \implies y \in \im$. Therefore
    \[P \cap Q = (\im/x\cdot \im)\cap (x\cdot \krn)/(x\cdot \im) = 0\]
    so the Lemma applies. The identifications
    \[M/P = (\krn/x\cdot \im)/(\im/x\cdot \im) \cong \krn/\im,\quad M/Q = (\krn/x\cdot \im)/(x\cdot\krn/x\cdot\im)\cong \krn/x\cdot \krn\]
    both follow from the third isomorphism theorem. Finally, we need to make an identification
    \[M/(P+Q) = (\krn/x\cdot \im)/((\im/x\cdot\im)+(x\cdot\krn/x\cdot\im)) \cong (\krn/\im)/((x\cdot \krn)/(x\cdot \im)).\]
    The bijection between submodules of $\krn/x\cdot \im$ and submodules of $\krn$ containing $x \cdot \im$ implies that 
    \[\im/x\cdot\im + x\cdot \krn/x\cdot\im = (\im + x\cdot \krn)/x\cdot \im,\]
    hence 
    \[M/(P + Q) \cong \krn/(\im + x\cdot \krn).\]
    We obtain a short exact sequence
    \[0\rightarrow (\im + x\cdot \krn)/\im \rightarrow \krn/\im \rightarrow \krn/(\im+ x\cdot \krn) \rightarrow 0.\]
    From here we make the identifications
    \[(\im + x\cdot \krn)/\im \cong (x\cdot\krn)/(\im \cap x\cdot \krn) = (x\cdot \krn)/(x\cdot \im).\]
    where the isomorphism is from the second isomorphism theorem. Hence we have the desired identification 
    \[M/(P+Q) \cong (\krn/\im)/(x\cdot\krn/x\cdot\im) \cong \krn/\im \otimes R/(x),\]
    since $x\cdot \krn/x\cdot \im = x \cdot(\krn/\im)$ again by $x\cdot \krn \cap \im = x \cdot \im.$
\end{proof}
\begin{remark}
    This yields a different way to Corollary 3.2 to see that if $f_1,\dots,f_n$ form a regular sequence that $\krn(\partial_{\underline{f}}^q)/x(\im(\partial_{\underline{f}}^{q-1}))$ is isomorphic to $\krn(\partial^q_{\underline{f}})\otimes R/(x)$, as the bottom row of the above square vanishes under that assumption.
\end{remark}
We want to extend these results to a global setting and make a geometric interpretation. In the proof of Proposition \ref{prop: affine}, we made use of the identity $\partial_{\underline{xf}} = x\partial_{\underline{f}}$. We make a similar identification in the global case.
 
\begin{lemma}
\label{decomp}
Suppose that we have a non-zero morphism $s: \mathcal{F}\rightarrow \Ocal_X$ where $\mathcal{F}$ is a locally free sheaf of constant rank. Let $\mathcal{L}$ be an invertible sheaf with section $\lambda: \Ocal_X \rightarrow \mathcal{L}$. Suppose furthermore that there is a map $t: \mathcal{F} \otimes \mathcal{L} \rightarrow \Ocal_X$ making the triangle commute
\begin{center}
\begin{equation}\label{eqn: triangle}
\begin{tikzcd}[row sep = 0.1em]\tag{$\ast$}
&& \mathcal{L}\\
&\mathcal{F}\otimes \mathcal{L}\arrow[ur,"s\otimes 1"]\arrow[dr,"t"]\\
&&\Ocal_X.\arrow[uu,"\lambda"]\\
\end{tikzcd}
\end{equation}
\end{center}
 Denote by $\partial$ the differential in $K^\bullet(\mathcal{F},s)$ and by $\delta$ the differential in $K^\bullet\left((\mathcal{F}\otimes \mathcal{L}), t\right) $. Then for any $n$,
 \[\krn(\partial^{-q} \otimes 1_{\mathcal{L}^n}) = \krn(\delta^{-q}) \otimes \mathcal{L}^{n-q},\quad \im(\partial^{-q-1}\otimes 1_{\mathcal{L}^n}) = \im(\delta^{-q-1}) \otimes \lambda(\Ocal_X) \otimes \mathcal{L}^{n-q+1}.\]
 \end{lemma}
 \begin{proof}
     We first show that for any $n$,
\[\begin{tikzcd}
    &\bigwedge^q(\mathcal{F}) \otimes \mathcal{L}^{n} \arrow[r,"\partial^{-q}\otimes 1"] &\bigwedge^{q-1}(\mathcal{F}) \otimes \mathcal{L}^{n}
\end{tikzcd}\]
can be expressed as the image of the morphism 
\[\begin{tikzcd}
 &\bigwedge^q (\mathcal{F}\otimes \mathcal{L}) \otimes \Ocal_X  \arrow[r,"\delta^{-q} \otimes \lambda"]
 &\bigwedge^{q-1}(\mathcal{F}\otimes \mathcal{L}) \otimes \mathcal{L}
\end{tikzcd}\]
\noindent
under the functor $(-\otimes \mathcal{L}^{n-q}).$ The commutativity of (\ref{eqn: triangle}) implies there is a commuting square
\[\begin{tikzcd}
    &(\mathcal{F}\otimes \mathcal{L})^q \arrow[d,"\sim"]\arrow[r,"1\otimes t"] &(\mathcal{F}\otimes \mathcal{L})^{q-1}\otimes \Ocal_X \arrow[r,"1\otimes \lambda"] &(\mathcal{F}\otimes \mathcal{L})^{q-1}\otimes \mathcal{L}\arrow[d,"\sim"]\\
    &\mathcal{F}^q\otimes \mathcal{L}^q \arrow[rr,"(1\otimes s)\otimes 1"] &&\mathcal{F}^{q-1}\otimes \mathcal{L}^q.
\end{tikzcd}\]
By definition of $1\wedge t, 1\wedge s$ and viewing $(\mathcal{F}\otimes \mathcal{L})^q$ as a subsheaf of $\bigwedge^q(\mathcal{F}\otimes \mathcal{L})$ etc (see the beginning of $\S2.2.1$), there is an induced commuting square of subsheaves
\[\begin{tikzcd}
    &\bigwedge^q(\mathcal{F}\otimes \mathcal{L}) \arrow[d,"\sim"]\arrow[r,"1\wedge t"] &\bigwedge^{q-1}(\mathcal{F}\otimes \mathcal{L})\otimes \Ocal_X \arrow[r,"1\otimes \lambda"] &\bigwedge^{q-1}(\mathcal{F}\otimes \mathcal{L})\otimes \mathcal{L}\arrow[d,"\sim"]\\
    &(\bigwedge^q\mathcal{F})\otimes \mathcal{L}^q \arrow[rr,"(1\wedge s)\otimes 1"] &&(\bigwedge^{q-1}\mathcal{F})\otimes \mathcal{L}^q.
\end{tikzcd}\]
The maps $1\wedge t$ and $1 \wedge s$ are the differentials in $K^\bullet(\mathcal{F}\otimes \mathcal{L},t)$ and $K^\bullet(\mathcal{F},s)$, respectively, and the induced vertical isomorphisms are from Lemma \ref{lem:line}. Hence applying the functor $(-\otimes \mathcal{L}^{n-q})$ yields the result. 
\par
Since $s \neq 0$ and a morphism of line bundles can either be injective or $0$, we have that $\lambda$ is injective by commutativity of (\ref{eqn: triangle}). As $\bigwedge^{q-1}\mathcal{F}\otimes \mathcal{L}$ is locally free, $1\otimes \lambda$ is injective. Hence 
    \[\krn(\delta^{-q}\otimes \lambda) = \krn((1\otimes \lambda)\circ(\delta^{-q}\otimes 1_{\Ocal_X})) = \krn(\delta^{-q}).\]
    As $\mathcal{L}$ is also locally free we therefore have 
    \[\krn(\partial^{-q}\otimes 1_{\mathcal{L}^n}) = \krn(\delta^{-q}\otimes \lambda \otimes 1_{\mathcal{L}^{n-q}}) = \krn(\delta^{-q}) \otimes \mathcal{L}^{n-q}.\]
    Similarly, $\im(\partial^{-q-1}\otimes 1_{\mathcal{L}^n}) = \im(\delta^{-q-1}\otimes \lambda \otimes 1_{\mathcal{L}^{n-q+1}}) = \im(\delta^{-q-1}) \otimes \lambda(\Ocal_X) \otimes \mathcal{L}^{n-q+1}$.
\end{proof}

We are now ready to prove the main result of this section. 
Let $X$ be a nonsingular variety and $Y_1,\dots,Y_n$ local complete intersection subvarieties of $X$. Assume that $Y = \bigcap Y_i$ decomposes into two lci components; $D$ of codimension 1 and  $Z$ of codimension $\geq 2$. Assume further that $Z = \bigcap Z_i$ where $Z_i = Y_i - D$ and that $D$ does not contain any irreducible components of $Z$. Let 
\begin{center}
\begin{equation}\label{eq: fibre}\tag{+}
    \begin{tikzcd}
        &Z \cap D \arrow[dr,"h"]\arrow[d,"i_Z"']\arrow[r,"j_D"] & D \arrow[d,"i"]\\
        &Z \arrow[r,"j"'] &X. 
    \end{tikzcd}
\end{equation}
    
\end{center}
 be the fibre square of closed immersions. We have the following 
 \begin{lemma}\label{Tor-ind}
    Let $X$ be a smooth variety and $D$ an effective Cartier divisor. Let $Z$ be a subvariety such that $D$ does not contain any associated points of $Z$. Then the fiber square (\ref{eq: fibre}) is Tor-independent. 
\end{lemma}
\begin{proof}
    Tor-independence is a local criterion, so we are reduced to proving the algebraic statement that 
    \[Tor^q_R(R/I,R/(x)) = 0 \quad \forall q\geq 1\]
    for $R$ a local ring and $x$ a non-zero divisor in $R$. There is a free resolution of $R/(x)$ given by 
    \[\begin{tikzcd}
        &R \ar[r,"x"] &R, 
    \end{tikzcd}\]
    so $Tor^2_R(R/I,R/(x)) = 0.$ For $q = 1$, we must show that the map 
    \[\begin{tikzcd}
        &R/I \ar[r,"x"] &R/I 
    \end{tikzcd}\]
    is injective. This follows from the fact that $D$ does not contain any associated points of $Z$ so that $x$ does not belong to the annihilator of any element of $R/I$. 
\end{proof}

\noindent
 We have global Koszul resolutions of each $\Ocal_{Y_i}$ and each $\Ocal_{Z_i};$
\[\mathcal{F}_i^\bullet = \{\Ocal(-Y_i) \xrightarrow{\sigma_{Y_i}} \Ocal_X \} ,\]
\[\mathcal{G}^\bullet_i = \{\Ocal(-Z_i) \xrightarrow{\sigma_{Z_i}} \Ocal_X \}.\]
  Consider the complex $\mathcal{G}^\bullet_1 \otimes \dots \otimes\mathcal{G}^\bullet_n$, which is a global Koszul complex whose cohomologies model $Tor^q_{\Ocal_X}(\Ocal_{Z_1},\dots,\Ocal_{Z_n})$. As $Z$ is lci, the excess sheaf $\mathcal{E}_Z$ is locally free and by \cite{https://doi.org/10.48550/arxiv.1510.04889} the cohomology sheaves of $\mathcal{G}_1^\bullet \otimes \dots \otimes \mathcal{G}_n^\bullet$ are $j_*\bigwedge^q \mathcal{E}_Z$. Denote its differentials by $\delta^q$. By Tor-independence we can identify $j_*(\bigwedge^q \mathcal{E}_Z|_{Z\cap D})$ with $i_*i^*j_*\bigwedge^q\mathcal{E}_Z$. Then we have the following 

\begin{theorem}
\label{main}
With notation as above
\[Tor_{\Ocal_X}^q(\Ocal_{Y_1},\dots,\Ocal_{Y_n}) = \mathcal{H}^q \otimes \Ocal(D)^{-q}\]
where $\mathcal{H}^q$ is the categorical pullback in the category $Mod_X$
\begin{center}
\begin{equation}\label{eq:square}\tag{$\dagger$}
     \begin{tikzcd}
&\mathcal{H}^q \arrow[r]\arrow[d] &i_*i^*\krn(\delta^{-q}) \arrow[d] \\
&j_*\bigwedge^q\mathcal{E}_Z \arrow[r] &j_*(\bigwedge^q\mathcal{E}_Z|_{Z\cap D}),
\end{tikzcd}
\end{equation}
   
\end{center}
where the lower horizontal arrow is $j_*$ applied to the restriction projection and the rightmost vertical arrow comes from $i_*i^*$ applied to the projection $\krn(\delta^{-q}) \rightarrow \krn(\delta^{-q})/\im(\delta^{-q-1})$.
\end{theorem}

\begin{proof}
Consider the complexes
\[\mathcal{F}_i^\bullet \otimes \Ocal(D) = \{\Ocal(-Z_i) \xrightarrow{\lambda \circ \sigma_{Z_i}} \Ocal(D)\}\]
where $\lambda : \Ocal_X \rightarrow \Ocal(D)$ is the defining section for $D$. Since each $\mathcal{F}_i^\bullet = K^\bullet(\Ocal(-Y_i),\sigma_{Y_i})$, we compute that 
\[\mathcal{C}^\bullet := \bigotimes (\mathcal{F}_i^\bullet \otimes \Ocal(D)) \cong K^\bullet\left(\bigoplus \Ocal(-Y_i),\sum \sigma_{Y_i}\right) \otimes \Ocal(D)^n.\]
The $q^{th}$ cohomology of $\mathcal{C}^\bullet$ is given by 
\[Tor^q_{\Ocal_X}(\Ocal_{Y_1},\dots,\Ocal_{Y_n}) \otimes \Ocal(D)^n\]
since $\Ocal(D)$ is locally free. 
\medskip
\newline
Now, in Lemma \ref{decomp}, take $\mathcal{F} = \bigoplus \Ocal(-Y_i), \mathcal{L} = \Ocal(D), \lambda = \text{ defining section of $D$}, s = \sum \sigma_{Y_i}, t = \sum \sigma_{Z_i}$. These data satisfy the of 
Lemma \ref{decomp}. Hence, denoting by $\partial$ the differentials in $\mathcal{C}^\bullet$ and by $\delta$ the differentials in $K^\bullet(\bigoplus \Ocal(-Z_i),\sum \sigma_{Z_i})$, we have
\[\krn(\partial^{-q}) = \krn(\delta^{-q}) \otimes \Ocal(D)^{n-q}, \quad \im(\partial^{-q-1}) = \im(\delta^{-q-1})\otimes \lambda(\Ocal_X) \otimes \Ocal(D)^{n-q-1}.\] 
\noindent
Then we have the identifications
\begin{align*}
    Tor^q_{\Ocal_X}(\Ocal_{Y_1},\dots,\Ocal_{Y_n}) \otimes \Ocal(D)^n \cong H^{-q}(\mathcal{C}^\bullet)= \frac{\krn(\partial^{-q})}{\im(\partial^{-q-1})}  &= \frac{\krn(\delta^{-q}) \otimes \Ocal(D)^{n-q}}{\im(\delta^{-q-1})\otimes \lambda(\Ocal_X) \otimes \Ocal(D)^{n-q-1}} \\
    &\cong \frac{\krn(\delta^{-q})}{\im(\delta^{-q-1})\otimes \Ocal(-D)} \otimes \Ocal(D)^{n-q}.
\end{align*}
\noindent
 Let $\mathcal{H}^q$ be the quotient $\frac{\krn(\delta^{-q})}{\im(\delta^{-q-1})\otimes \Ocal(-D)}.$

\medskip
To see $\mathcal{H}^q$ as the pullback in the diagram (\ref{eq:square}), we first construct a morphism from $\mathcal{H}^q$ to the pullback and then check that it is locally an isomorphism. We note that factoring out the inclusions of $\krn(\delta^{-q})\otimes \Ocal(-D)$ and $\im(\delta^{-q-1})$ into $\krn(\delta^{-q})$ induces projections
\[\mathcal{H}^q \rightarrow \krn(\delta^{-q}) \otimes \Ocal_D \cong i_*i^* \krn(\delta^{-q}), \quad \mathcal{H}^q \rightarrow \krn(\delta^{-q})/\im(\delta^{-q-1}) \cong j_*\bigwedge^q\mathcal{E}_Z.\]
These projections are the morphisms we use in the square (\ref{eq:square}). With these morphisms, the square (\ref{eq:square}) locally becomes the commutative square in Proposition \ref{alg-pullback}, so (\ref{eq:square}) commutes globally. We therefore have a global induced map from $\mathcal{H}^q$ to the pullback. Since locally the square (\ref{eq:square}) becomes the square in \ref{alg-pullback}, we need to check the condition that $x$ does not belong to any prime ideal in $Ass_R(\krn(\delta^{-q})/\im(\delta^{-q-1}))$, where $x$ is the local defining section for $D$. As 
$Tor^q_{\Ocal_X}(\Ocal_{Z_1},\dots,\Ocal_{Z_n}) = \krn(\delta^{-q})/\im(\delta^{-q-1})$ is locally free on $Z$, its associated points are just the associated points of $\Ocal_Z$, which by Lemma \ref{no-embedded} correspond to the generic points of the irreducible components of $Z$. Since $D$ does not contain any of these points by assumption, its local defining section $x$ cannot belong to any of the associated primes of the cohomology modules. Therefore, by Proposition \ref{alg-pullback}, the square (\ref{eq:square}) is locally a cartesian square. Thus our global morphism from $\mathcal{H}^q$ to the pullback is locally an isomorphism.
\end{proof}

\begin{example}\label{Example 2}
We give an example to show that our answer is not in the form of exterior powers. Consider a line through a hyperplane in 4-space, we will use the example of the line $\{y=z=w = 0\}$ and the hyperplane $\{x = 0\}$ in $\mathbb{A}^4_k$, considered as the intersection of the hypersurfaces $\{xy = 0\}, \{xz = 0\}, \{xw = 0\}$. Let $R = k[x,y,z,w]$. We are finding the cohomologies of the Koszul complex $K^\bullet(xy,xz,xw)$. By \ref{prop: affine} we have 
\[H^{-1} \cong (R/(x))^3/((w,-z,y)),\quad H^{-2} \cong R/(x).\]
We compute the second exterior power of $H^{-1}$. Pick generators $e_1,e_2,e_3$ of $H^{-1}$ subject to the relation $ze_2 = we_1 + ye_3$. Then $\bigwedge^2H^{-1}$ is generated by $e_1\wedge e_2,e_1\wedge e_3, e_2\wedge e_3$ which are subject to the relations
\[z(e_1\wedge e_2) = y(e_1 \wedge e_3), z(e_2\wedge e_3) = w(e_1 \wedge e_3), w(e_1 \wedge e_2) = y(e_2 \wedge e_3).\]
Hence 
\[\bigwedge^2 H^{-1} \cong (R/(x))^3/((-z,y,0), (-w,0,y), (0,-w,z)) \cong (y,z,w) \subset R/(x).\]
Note that this shows that $\bigwedge^2H^{-1}$ is a free $R/(x)$-module of rank one \textit{away from the line \{y=z=w=0\}}, and therefore coincides with $H^{-2}$ there. So, as expected, the multitors do form an exterior algebra \textit{away from the intersection of the components.}
\end{example}




\section{Further Work}
We wish to make some remarks about generalisations of the previous results. Let $X$ be a non-singular variety and $Y_1,\dots,Y_l$ local complete intersection subvarieties whose intersection is $\bigcap Y_i = \bigcup Z_j$ can be decomposed into a union of local complete intersection components (not necessarily irreducible). Pick any $Z_j$. Then on the set $U_j = X\backslash \bigcup Z_k$ where the union is taken over all $k \neq j$ we are in the situation of Scala's Theorem. Thus 
\[Tor^q_{\Ocal_X}(\Ocal_{Y_1},\dots,\Ocal_{Y_l})|_{U_j} \cong i_{Z_j*}\bigwedge^q \mathcal{E}_{Z_j}|_{U_j}\]
where $i_{Z_j}$ is the inclusion of $Z_j$ in $X$. Hence, whatever our result will be globally, it will need to restrict to each component as the exterior powers of the excess bundle for that component, and the intersections of the components are the only places where there is interesting behaviour. However, there is an imprecision with this statement, which is that the component-wise excess sheaves $\mathcal{E}_{Z_j}$ are only defined on the complement of the other components. While there are cases where one may extend these excess sheaves to locally free sheaves on the whole of $Z_j$, even in our case this will not yield the correct answer. Indeed, if there is no excess on $Z$, our theorem gives the answer as $\krn(\delta^{-q})\otimes \Ocal_D$ which is not necessarily a locally free $\Ocal_D$-module. Given the nature of our answer, it seems natural to believe that there will be projections from the global multitor to \textit{coherent} extensions of the excess bundles of each $Z_i$. We conjecture that the multitors will form limits over systems consisting of these projections and further projections to coherent sheaves supported on the intersections of the lci components $Z_j$. 


\nocite{*}
\printbibliography


\end{document}